\documentclass[journal]{IEEEtran}
\IEEEoverridecommandlockouts
\usepackage{amsmath,amssymb,bm,mathtools}
\usepackage{graphicx,booktabs,cite,microtype,url}
\usepackage[hidelinks]{hyperref}
\hypersetup{pdftitle={Operational Interpretation of Controllability Allocation Scores via Randomized Interventions},pdfauthor={Kazuhiro Sato}}
\newtheorem{theorem}{Theorem}
\newtheorem{proposition}{Proposition}
\newtheorem{corollary}{Corollary}
\newtheorem{remark}{Remark}
\newcommand{\R}{\mathbb R}
\newcommand{\E}{\mathbb E}
\newcommand{\Pset}{\mathcal P}
\newcommand{\Dset}{\mathcal D}

\newcommand{\norm}[1]{\lVert#1\rVert}
\newcommand{\Gbar}{\overline G}

\newcommand{\argmin}{\mathop{\rm argmin}\limits}

\mathtoolsset{showonlyrefs=true}

\DeclareMathOperator{\tr}{tr}
\DeclareMathOperator{\diag}{diag}

\begin{document}
\title{Operational Interpretation of Controllability Allocation Scores via Randomized Interventions}
\author{Kazuhiro Sato%
\thanks{K. Sato is with the Department of Mathematical Informatics, Graduate School of Information Science and Technology, The University of Tokyo, Tokyo 113-8656, Japan (e-mail: kazuhiro@mist.i.u-tokyo.ac.jp).}}
\maketitle

\begin{abstract}
Controllability allocation scores (CASs) assign weights to
virtual input channels without specifying how interventions
should be implemented.
This paper provides an operational interpretation through
randomized finite-duration interventions with prescribed,
possibly candidate-dependent waveforms.
Selection probabilities are fixed before a mean-output
target is specified, while waveform amplitudes are
optimized after the target is revealed.
We derive the unique minimum expected-energy amplitude
policy and formulate volumetric and weighted average-energy
allocation criteria.
We also establish a finite-trial mean-output error bound
and show that waveform-based weighted AE-CAS minimizes
an upper bound on the target-averaged mean-square error.
We establish
uniform convergence of normalized waveform-response
Gramians and show that accumulation points of optimal
allocations solve the limiting problem.
With uniform timing, equal positive waveform efficiencies
recover the original CAS in the short-duration limit,
whereas unequal efficiencies yield waveform-weighted limits.
Under the former conditions, the limiting error bound
is proportional to the original weighted AE-CAS objective.
For common waveforms, we further establish an exact
finite-duration correspondence with the original
full-state V-CAS under uniform timing.
\end{abstract}

\section{Introduction}
\label{sec:intro}

Controllability allocation scores (CASs) provide
an optimization-based framework for allocating
limited control capability among prescribed
input directions to influence designated
target variables in dynamical
networks~\cite{SatoCAS2026}.
Building on energy-based controllability
metrics~\cite{baggio2022energy}, CAS accounts
for how network dynamics transmit the effects
of candidate inputs to the desired outputs.
Specifically, CAS considers the system
\begin{equation}
\label{eq:cas-virtual-system}
\begin{aligned}
\dot x(t)&=Ax(t)+B(p)u(t),\quad y(t)=Cx(t),\\
B(p)&=B_0\operatorname{diag}
(\sqrt{p_1},\ldots,\sqrt{p_m}),
\end{aligned}
\end{equation}
where $A\in\R^{n\times n}$,
$B_0=[b_1\ \cdots\ b_m]\in\R^{n\times m}$
specifies the candidate input directions,
and $C\in\R^{r\times n}$ has full row rank
and specifies the target outputs.
The allocation $p$ belongs to a nonempty
compact convex set
$\Pset\subseteq\Delta_m
:=\{p=(p_i)\in\R^m:p\geq0,\ p_1+ \cdots +p_m=1\}$,
where $m$ and $r$ need not coincide.
For a horizon $T>0$, the output controllability
Gramian is
\begin{equation}
\label{eq:cas-gramian}
\begin{aligned}
W(p,T)&:=\sum_{i=1}^m p_iW_i(T),\\
W_i(T)&:=\int_0^T
Ce^{At}b_i b_i^\top e^{A^\top t}C^\top\,{\rm d}t.
\end{aligned}
\end{equation}
The volumetric CAS (V-CAS) and average-energy
CAS (AE-CAS) are defined, respectively, by
\begin{equation}
\label{eq:cas-objectives}
\begin{aligned}
p_{\mathrm V}^\star
&\in\argmin_{p\in\Pset_T}
-\log\det W(p,T),\\
p_{\mathrm A}^\star
&\in\argmin_{p\in\Pset_T}
\operatorname{tr}\bigl(W(p,T)^{-1}\bigr),
\end{aligned}
\end{equation}
where $\Pset_T:=\{p\in\Pset:W(p,T)\succ0\}$
is assumed nonempty.
The V-CAS maximizes the volume of reachable
output targets under a fixed control-energy
budget, whereas the AE-CAS minimizes the
average minimum energy required to reach
unit-norm output targets.
The conventional controllability scores~\cite{SatoTerasaki2024,SatoKawamura2025}
are recovered when $B_0=C=I_n$ and
$\Pset=\Delta_n$.

However, in the CAS formulation, the allocation weights $p_1, \ldots, p_m$ determine the gains of virtual input channels in system \eqref{eq:cas-virtual-system} without specifying how interventions should be implemented in practice. Repeated neural perturbation experiments provide a motivating setting. Optogenetic stimulation can target selected neurons or neuronal ensembles across trials \cite{triplett2025rapid}, while trial-averaged responses can reveal stimulus-dependent population activity patterns despite variability across individual trials \cite{mulholland2024self}. Moreover, optogenetic control experiments have demonstrated the feasibility of eliciting prescribed neural firing-rate patterns \cite{bolus2018design}. Since the significance of trial-averaged responses depends on the experimental context \cite{tlaie2024does}, we consider settings in which the mean response is itself an explicit control objective.

\begin{table*}[t]
\centering
\caption{Two-stage randomized intervention framework}
\label{tab:two-stage}
\renewcommand{\arraystretch}{1.35}
\begin{tabular*}{\textwidth}{
@{\extracolsep{\fill}}
p{0.12\textwidth}
p{0.20\textwidth}
p{0.27\textwidth}
p{0.30\textwidth}
@{}}
\hline
Stage & Information & Decision & Implementation \\
\hline
\textbf{1: Allocation}
& Specific target $\eta$ unknown
& Choose $p\in\Pset$.
& Fix the selection probabilities
  $\Pr\{I=i\}=p_i$. \\
\hline
\textbf{2: Intervention}
& Target $\eta$ revealed;
  allocation $p$ fixed
& Choose the amplitude policy
  $u_i(s;\eta,p)$ to minimize
  expected input-signal energy.
& In each trial, sample $(I,S)$
  and apply one prescribed waveform, subject to
  $\mathbb E_p[D]=\eta$. \\
\hline
\end{tabular*}
\end{table*}

These experimental settings motivate the question of how CAS weights can be translated into a stimulation strategy that accommodates different future mean-output targets. Previous studies have investigated Gramian-based actuator selection \cite{summers2016}, randomized actuator scheduling \cite{siami2021}, and model-based design of neural stimulation sites and inputs \cite{de2026}. While these studies address the selection, timing, and design of control inputs, they do not directly establish how CAS weights can be interpreted as selection probabilities across repeated intervention trials. To address this question, we investigate a randomized finite-duration intervention framework in which a common allocation of intervention opportunities is used for different future mean-output targets.

The contributions of this paper are summarized as follows.
\begin{enumerate}
\item We provide an operational interpretation of CAS
through a two-stage randomized intervention framework
with prescribed finite-duration waveforms.
The allocation determines the probabilities of selecting
candidate input directions before a mean-output target
is known, while waveform amplitudes are optimized
after the target is revealed.
We derive the unique minimum expected-energy
amplitude policy and the resulting volumetric and
weighted average-energy allocation criteria.

\item We establish theoretical relationships between
waveform-based and original CAS.
We characterize the short-duration limit of the
waveform-response Gramian and establish convergence
properties of optimal allocations.
We introduce waveform efficiency and show that,
under uniform timing, equal positive efficiencies
recover the original CAS, whereas unequal
efficiencies yield waveform-weighted limits.
For common waveforms, we further establish
an exact finite-duration correspondence
with the original full-state V-CAS.
These results suggest that the original CAS can guide preliminary intervention-probability design before stimulation waveforms are finalized, potentially helping prioritize neural stimulation sites for subsequent calibration. Such a staged approach is related to the identification of stimulation sites before stimulation-pattern design discussed in \cite{de2026}.

\item  We establish a finite-trial mean-output error
bound and show that waveform-based weighted AE-CAS
minimizes an upper bound on the target-averaged
mean-square error. Under uniform timing and equal
positive efficiencies, the limiting upper bound
is proportional to the original weighted
AE-CAS objective.
\end{enumerate}

The rest of this paper is organized as follows.
Section~II introduces the randomized intervention
framework and derives the waveform-based allocation
criteria. Section~III establishes the short-duration
limits, waveform-dependent CAS recovery, and an exact finite-duration correspondence.
Section~IV presents numerical examples,
and Section~V concludes the paper.

\emph{Notation:}
The symbols $(\cdot)^\top$ and $\tr$ denote
transpose and trace, respectively.
The norm $\norm{\cdot}$ denotes the Euclidean norm
for vectors and the induced spectral norm for matrices,
while $\norm{\cdot}_{\rm F}$ denotes the Frobenius norm.
For a symmetric matrix $X$, $\lambda_{\min}(X)$
denotes its smallest eigenvalue.
The notation $X\succ0$ means
that $X$ is positive definite.

\section{Randomized Interventions}
\label{Sec2}
\label{sec:scores}

The proposed framework fixes the probabilities
of selecting input directions before the
mean-output target $\eta$ is known (Stage 1)
and optimizes waveform amplitudes after
$\eta$ is revealed (Stage 2), as summarized
in Table~\ref{tab:two-stage}.
We first solve the minimum expected-energy
problem in Stage 2 to formulate the
Stage 1 allocation criteria.
We then establish a finite-trial mean-output
error bound.

\subsection{Randomized Intervention Protocol}

We retain $A$, $B_0$, $C$, and $\Pset$ from Section~I, but now
interpret $p_i$ as the probability of selecting direction $b_i$ in
an intervention trial, rather than as a virtual input weight.
Fix a terminal time $T>0$ and a intervention duration
$0<\delta<T$. For each candidate $b_i$, prescribe a real waveform
$\psi_i\in L^2([0,1])$, extended by zero outside $[0,1]$, with
\begin{align}
\int_0^1\psi_i(v)^2\, {\rm d}v=1. \label{psi_assumption}    
\end{align}
The shapes of $\psi_i$ are fixed independently of $p$, $\eta$, and $\delta$.
This normalization entails no loss for nonzero prescribed
shapes, because their norms can be absorbed into the unrestricted
amplitudes. 
For a given allocation $p\in\Pset$, each trial
samples a candidate direction $I$ and a intervention
starting time $S$ according to the joint
distribution $\mu_p$ defined by
\begin{equation}
 \mu_p(\{i\}\times B)
 =p_i\int_B q_{i,\delta}(s)\,{\rm d}s,
 \label{eq:law}
\end{equation}
for every $i\in\{1,\ldots,m\}$ and
measurable set $B\subseteq[0,T-\delta]$,
where $q_{i,\delta}$ is a probability density
on $[0,T-\delta]$ that may depend on $\delta$
but is fixed independently of $p$ and $\eta$.
The interval $[0,T-\delta]$ ensures that every
intervention ends by $T$.

The Stage 2 intervention protocol adapts the waveform amplitudes
to a revealed mean-output target while keeping the selection
law~\eqref{eq:law} fixed.
Specifically, after a target $\eta\in\R^r$ is revealed,
an amplitude policy $u_i(s;\eta,p)$ assigns a possibly signed
amplitude to each candidate direction $i$ and starting time $s$.
In a trial with $(I,S)=(i,s)$, the system evolves according to
\begin{equation}
  \dot x(t)=Ax(t)+b_i u_i(s;\eta,p)g_{i,\delta}(t;s),\quad
  y(t)=Cx(t),
 \label{eq:plant}
\end{equation}
where $g_{i,\delta}(t;s):=\psi_i((t-s)/\delta)$.
Thus, each trial applies one prescribed waveform along one
candidate direction.
Although all trials share the same system matrices $A$ and $C$ and candidate input directions $b_1,\ldots,b_m$, their initial states may differ.
We therefore evaluate intervention effects through terminal
output changes relative to the corresponding unforced trajectories.
For a fixed target and policy, trials are assumed independent,
and the target is required to be achieved in expectation rather
than in every individual trial.

\begin{remark}
The proposed framework can be interpreted as a simplified
model of interventions in neuronal populations.
Optogenetic techniques allow stimulation of different
neuronal ensembles across trials \cite{triplett2025rapid}.
Here, $b_i$ represents a candidate stimulation pattern
and $p_i$ its selection probability.
The target $\eta$ represents a desired trial-averaged
change in terminal neuronal population activity relative
to the corresponding unforced dynamics.
The model focuses on the allocation of stimulation
opportunities rather than the full complexity of
neuronal dynamics.
\end{remark}

\subsection{Mean-Output Steering and Minimum Energy}

We first characterize the terminal output
change produced by a single intervention.
For a trial with $(I,S)=(i,s)$ and
initial state $x_0$, \eqref{eq:plant} yields
\begin{align}
D&:=y(T)-Ce^{AT}x_0
  =k_i^\delta(s)u_i(s), \label{eq:response}
\\
k_i^\delta(s)
&:=\delta\int_0^1 h_i(s+\delta v)\psi_i(v)\,{\rm d}v,\quad
h_i(t) :=Ce^{A(T-t)}b_i, \nonumber
\end{align}
where the arguments $(\eta,p)$ of $u_i$ are suppressed.
The random pair $(I,S)$ induces variability in both $D$ and the
input-signal energy, so we evaluate policies in expectation.

For a fixed allocation $p$, the Stage 2 problem
seeks an amplitude policy that achieves the
prescribed mean-output change with minimum
expected input-signal energy.
To formulate this problem, let $\E_p$ denote
expectation with respect to $\mu_p$ defined by \eqref{eq:law}, and consider
square-integrable amplitude policies
$u\in L^2(\mu_p)$, where $u(i,s)=u_i(s)$.
The expected output change and input-signal
energy per trial are
\begin{align}
 \mathcal L_pu
 &:=\E_p[D]
 =\E_p[k_I^\delta(S)u_I(S)],
 \label{eq:mean}\\
 J_p(u)
 &:=\E_p\left[\int_0^T
 \left(u_I(S)g_{I,\delta}(t;S)\right)^2
 \,{\rm d}t\right]\\
 &=\delta\E_p[u_I(S)^2].
 \label{eq:energy}
\end{align}
The last equality follows from the waveform normalization \eqref{psi_assumption}.

For any $p$ with $G_\delta(p)\succ 0$, the minimum expected
input-signal energy required to achieve a target
$\eta\in\R^r$ is
\begin{equation}
 \min_{\mathcal L_pu=\eta}J_p(u)
 =\delta\eta^\top G_\delta(p)^{-1}\eta,
 \label{eq:value}
\end{equation}
where
$G_\delta(p)$ denotes the waveform-response
Gramian defined by
\begin{equation}
 \begin{aligned}
 G_\delta(p)
  &:=\sum_{i=1}^m p_iG_{\delta,i}
  =\E_p[k_I^\delta(S)k_I^\delta(S)^\top],\\
  G_{\delta,i}
  &:=\int_0^{T-\delta}q_{i,\delta}(s)
  k_i^\delta(s)k_i^\delta(s)^\top\,{\rm d}s.
\end{aligned}
 \label{eq:gram}
\end{equation}
The Gramian $G_\delta(p)$ represents the
second moment of the random waveform response
and depends affinely on $p$.
In this paper, we assume that
\begin{equation}
 \Dset_\delta
 :=\{p\in\Pset:G_\delta(p)\succ0\}
 \ne\varnothing.
 \label{eq:domain}
\end{equation}

To prove \eqref{eq:value},
we adapt the Hilbert-space minimum-norm argument
of Luenberger~\cite[Sec.~6.10]{Luenberger1969}
to the present randomized intervention setting,
where amplitude policies belong to $L^2(\mu_p)$
and the constraint is imposed on the expected
terminal-output change.
For any fixed $p\in\Dset_\delta$ and
$\eta\in\R^r$, the minimum in
\eqref{eq:value} is uniquely attained,
up to $\mu_p$-almost-everywhere equality,
by the amplitude policy
\begin{equation}
 u_i^*(s;\eta,p)
 =k_i^\delta(s)^\top G_\delta(p)^{-1}\eta.
 \label{eq:policy}
\end{equation}
To verify this, the definition \eqref{eq:gram} of
$G_\delta(p)$ shows that
$\mathcal L_pu^*=\eta$.
For any feasible policy $u\in L^2(\mu_p)$,
we have
$\E_p[k_I^\delta(S)(u_I(S)-u_I^*(S))]=0$,
which, together with \eqref{eq:policy}, yields
$\E_p[u_I^*(S)(u_I(S)-u_I^*(S))]=0$.
Consequently,
\[
J_p(u)
=J_p(u^*)
+\delta\E_p[(u_I(S)-u_I^*(S))^2]
\ge J_p(u^*).
\]
Since $\delta>0$, equality holds if and only if
$u=u^*$ $\mu_p$-almost everywhere.
Thus, substituting \eqref{eq:policy} into \eqref{eq:energy},
we obtain
\eqref{eq:value}.

\subsection{Waveform-Based Allocation Criteria}

The minimum-energy expression~\eqref{eq:value} determines how to
allocate interventions before future mean-output targets $\eta$ are specified.
We consider two criteria according to the available target information.
\begin{itemize}

\item \emph{No target distribution specified.}
When no probability distribution of future
targets $\eta$ is specified, we consider
the range of targets attainable under a fixed
expected input-signal energy budget $\rho^2>0$.
For $p\in\Dset_\delta$, these targets form
the ellipsoid
 $\{\eta\in \mathbb{R}^r:\delta\eta^\top G_\delta(p)^{-1}\eta
 \leq\rho^2\}$, 
whose volume is proportional to
$\sqrt{\det G_\delta(p)}$.
Thus, maximizing the volume is equivalent
to maximizing $\det G_\delta(p)$.
This distribution-free criterion may be useful
in exploratory neural perturbation studies,
where the behavioral effects of stimulating neuronal ensembles
are investigated experimentally~\cite{bounds2025network,carrillo2019controlling}.
It maximizes the volume of feasible mean targets, not the efficiency
of discovering behaviorally useful patterns.

\item \emph{Statistical information available.}
When the second-moment matrix
$Q:=\E_\eta[\eta\eta^\top]\succ 0$
of future targets $\eta$ is known, the allocation
can be chosen to minimize the average
minimum expected input-signal energy
over future targets.
Averaging~\eqref{eq:value} over future targets
yields
\begin{align}
\E_\eta\left[\min_{\mathcal L_pu=\eta}J_p(u)\right]
 =\delta\tr(QG_\delta(p)^{-1}). \label{Energy_AECAS}
\end{align}
This criterion accounts for both the mean and
variability of future targets through $Q$,
without requiring their full distribution,
detailed
in \cite[Sec.~II-C]{sato2026task}.

\end{itemize}

Both criteria can be expressed in the
common form
\begin{equation}
 p_{\nu,\delta}
 \in\argmin_{p\in\Dset_\delta}
 \phi_\nu(G_\delta(p)),
 \qquad \nu\in\{V,A\},
 \label{eq:scoring}
\end{equation}
where, for $X\succ0$,
\[
 \phi_V(X):=-\log\det X,\qquad
 \phi_A(X):=\tr(QX^{-1}).
\]
The volumetric criterion $\phi_V$ defines the
waveform-based V-CAS, while the average-energy
criterion $\phi_A$ defines a weighted extension of
the waveform-based AE-CAS, reducing to the
standard AE-CAS when $Q=cI_r$ for $c>0$ \cite{SatoCAS2026}.

\begin{remark}
If each trial selects $I=i$ with probability $p_i$ and replaces
$u_i(s;\eta,p)g_{i,\delta}(t;s)$ in~\eqref{eq:plant} by an arbitrary
$u_i\in L^2([0,T])$, then
$D_i=\int_0^T h_i(t)u_i(t)\,{\rm d}t$.
For $W(p,T)\succ0$, minimizing
$\E_{I\sim p}[\int_0^T u_I(t)^2\,{\rm d}t]$ subject to
$\E_{I\sim p}[D_I]=\eta$ gives
$u_i^*(t)=h_i(t)^\top W(p,T)^{-1}\eta$ and minimum value
$\eta^\top W(p,T)^{-1}\eta$, by the same orthogonal-decomposition
argument as above. Thus, at every finite horizon, the Stage~1 problems
have exactly the same optimal allocations as the original V-CAS and
the corresponding weighted AE-CAS.
\end{remark}

\subsection{Finite-Trial Mean-Output Accuracy} \label{subsec2-D}

The constraint $\mathcal L_pu=\eta$ ensures
that the expected terminal output change equals
$\eta$, but the empirical mean over finitely many
trials generally differs from this target.
The following theorem quantifies this deviation
and establishes its connection to the
waveform-based weighted AE-CAS.

\begin{theorem}
\label{thm:sampling}
Fix an allocation $p\in\Pset$, a target
$\eta\in\R^r$, and an amplitude policy
$u\in L^2(\mu_p)$ satisfying
$\mathcal L_pu=\eta$.
Let $D_1,\ldots,D_N$ be independent copies
of the terminal output change $D$ under
this allocation and policy. Define
\begin{align*}
\bar D_N:=\frac{\sum_{j=1}^N D_j}{N},\quad K_\delta:=\max_{i\in\{1,\ldots, m\},\,s\in[0,T-\delta]}
\norm{k_i^\delta(s)}.
\end{align*}
Then,
\begin{equation}
\E_p\norm{\bar D_N-\eta}^2
\leq\frac{K_\delta^2}{N\delta}J_p(u).
\label{eq:samplingerror}
\end{equation}
Furthermore, suppose that $p\in\mathcal D_\delta$
and the Stage 2 optimal policy \eqref{eq:policy} is used
for each future target $\eta$.
Then
\begin{equation}
\E_\eta\left[
\E_p\norm{\bar D_N-\eta}^2
\right]
\leq\frac{K_\delta^2}{N}
\tr(QG_\delta(p)^{-1}).
\label{eq:samplingaverage}
\end{equation}
\end{theorem}

\begin{IEEEproof}
Since $\mathcal L_pu=\eta$, we have
$\E_p[D_j]=\eta$ for every $j$.
By independence, for $j\neq k$,
$\E_p[(D_j-\eta)^\top(D_k-\eta)]=
\E_p[D_j-\eta]^\top
\E_p[D_k-\eta]=0$.
Consequently,
\[
\begin{aligned}
\E_p\norm{\bar D_N-\eta}^2
=\frac{1}{N^2}
\sum_{j=1}^N\E_p\norm{D_j-\eta}^2=\frac{1}{N}\E_p\norm{D-\eta}^2.
\end{aligned}
\]

Next, using $\E_p[D]=\eta$ and
\eqref{eq:response}, we obtain
\[
\begin{aligned}
\E_p\norm{D-\eta}^2
&\leq\E_p\norm{D}^2=\E_p[\norm{k_I^\delta(S)}^2u_I(S)^2]\\
&\leq K_\delta^2\E_p[u_I(S)^2]=\frac{K_\delta^2}{\delta}J_p(u).
\end{aligned}
\]
Combining these estimates proves
\eqref{eq:samplingerror}.
Furthermore, \eqref{eq:policy}, \eqref{Energy_AECAS}, and \eqref{eq:samplingerror} yield \eqref{eq:samplingaverage}.
\end{IEEEproof}

Theorem~\ref{thm:sampling} means that
the waveform-based weighted
AE-CAS minimizes the upper bound in
\eqref{eq:samplingaverage}.

\section{Short-Duration Limits and CAS Recovery}
\label{sec:limit}

The waveform-based allocation problem~\eqref{eq:scoring} retains the convex structure
of CAS problem~\eqref{eq:cas-objectives}, but their candidate Gramians~\eqref{eq:gram} depend on waveform shape,
duration, and intervention timing. We first establish the leading
short-duration limit and convergence of its optimal allocations.
We then identify conditions for recovering the original CAS.

\subsection{General Starting-Time Distributions}

The unit normalization \eqref{psi_assumption}
determines the scaling of the waveform-response Gramian as
$\delta\downarrow0$.
Indeed, with
\begin{align}
    M_0:=\max_{i\in\{1,\ldots, m\},t\in[0,T]}\norm{h_i(t)}, \label{def_M0}
\end{align}
Cauchy--Schwarz inequality and~\eqref{eq:response} give
$\norm{k_i^\delta(s)}
\leq\delta M_0\int_0^1|\psi_i(v)|\,{\rm d}v
\leq\delta M_0$.
Since $q_{i,\delta}$ is a probability density and $p\in\Pset$,
this and \eqref{eq:gram} yield
$\norm{G_\delta(p)}\leq\delta^2 M_0^2$
uniformly over $p\in\Pset$.
Thus, $G_\delta(p)$ converges to zero for every allocation
as $\delta\downarrow0$.

To identify the limiting dependence on $p$,
we remove this common vanishing factor by introducing
the normalized response kernel
$\bar k_i^\delta(s):=k_i^\delta(s)/\delta$
and the normalized Gramian
\begin{equation}
\begin{aligned}
 \Gbar_\delta(p)
 &:=\frac{G_\delta(p)}{\delta^2}
 =\sum_{i=1}^m p_i\Gbar_{\delta,i},\\
 \Gbar_{\delta,i}
 &:=\int_0^{T-\delta}
 q_{i,\delta}(s)
 \bar k_i^\delta(s)\bar k_i^\delta(s)^\top\,{\rm d}s.
\end{aligned}
\label{eq:normalizedwaveform}
\end{equation}
By \eqref{eq:response} and continuity of $h_i$, as $\delta\downarrow0$,
\[
\bar k_i^\delta(s)
\longrightarrow m_i h_i(s),
\qquad
m_i:=\int_0^1\psi_i(v)\,{\rm d}v,
\]
uniformly over $[0,T-\delta]$.
Consequently,
$\bar k_i^\delta(s)\bar k_i^\delta(s)^\top
\longrightarrow
m_i^2h_i(s)h_i(s)^\top$
uniformly over $[0,T-\delta]$.
Since $q_{i,\delta}$ is a probability density,
this and \eqref{eq:normalizedwaveform} yield
\begin{align}
\Gbar_{\delta,i}
&=\gamma_i
\int_0^{T-\delta}
q_{i,\delta}(s)h_i(s)h_i(s)^\top\,{\rm d}s
+o(1),\\
\gamma_i&:=m_i^2. \label{efficiency}
\end{align}
Thus, $\gamma_i$ quantifies the leading-order
contribution of the prescribed waveform to the
normalized Gramian, and we call it the
waveform efficiency.
The unit normalization \eqref{psi_assumption} and
Cauchy--Schwarz inequality imply
$0\leq\gamma_i\leq1$.

To characterize the limit as $\delta\downarrow0$,
we assume that the starting-time densities
$q_{i,\delta}$ admit limiting probability
densities on $[0,T]$ in the $L^1$ sense.
Specifically, let $\tilde q_{i,\delta}$ denote
the extension of $q_{i,\delta}$ by zero
on $(T-\delta,T]$.
We assume that there exist probability densities
$q_{i,0}$ on $[0,T]$ such that
\begin{equation}
\varepsilon_\delta:=
\max_{1\le i\le m}
\int_0^T
|\tilde q_{i,\delta}(s)-q_{i,0}(s)|\,{\rm d}s
\longrightarrow 0
\quad\text{as }\delta\downarrow0.
\label{eq:densityerror}
\end{equation}
We then define the waveform-efficiency-weighted
Gramian by
\begin{equation}
 \Gbar_0(p):=\sum_{i=1}^m p_i\Gbar_{0,i},\quad
 \Gbar_{0,i}:=\gamma_i\int_0^T q_{i,0}(s)
 h_i(s)h_i(s)^\top\,{\rm d}s.
\label{eq:weightedgramian}
\end{equation}

The following theorem establishes convergence
of optimal solutions to
\begin{equation}
 p_{\nu,\delta}
 \in\argmin_{p\in\Dset_\delta}
 \phi_\nu(\Gbar_\delta(p)),
 \qquad \nu\in\{V,A\},
 \label{eq:scoring2}
\end{equation}
where
$\Dset_\delta$ is defined by \eqref{eq:domain} and
$\Dset_0:=\{p\in\Pset:\Gbar_0(p)\succ0\}$.
 
\begin{theorem}
\label{thm:limit}
Suppose~\eqref{eq:densityerror} holds. With $M_0$ defined in \eqref{def_M0},
let
$\displaystyle M_1:=\max_{i\in\{1,\ldots, m\},\, t\in[0,T]}\norm{h_i'(t)}$.
Set $\displaystyle \beta_\psi:=\max_{i\in\{1,\ldots,m\}}\int_0^1v|\psi_i(v)|\,{\rm d}v$.
Then, for every $0<\delta<T$,
\begin{equation}
 \sup_{p\in\Delta_m}
 \norm{\Gbar_\delta(p)-\Gbar_0(p)}
 \leq 2\beta_\psi M_0M_1\delta+M_0^2\varepsilon_\delta.
 \label{eq:uniform}
\end{equation}
Suppose, in addition, that $p^\circ\in\Pset$
satisfies $\Gbar_0(p^\circ)\succ0$.
Then, for each $\nu\in\{V,A\}$,
\eqref{eq:scoring2} admits an optimizer
for all sufficiently small $\delta>0$.
Every accumulation point of optimal solutions
$p_{\nu,\delta}$ as $\delta\downarrow0$
minimizes $\phi_\nu(\Gbar_0(p))$ over $\Dset_0$.
If the limiting optimizer is unique,
every choice of optimal solutions $p_{\nu,\delta}$
converges to it as $\delta\downarrow0$.
\end{theorem}
\begin{IEEEproof}
To establish the uniform bound \eqref{eq:uniform}, we first estimate the difference between the normalized response kernel $\bar k_i^\delta(s)$ and its short-duration limit $m_i h_i(s)$.
The bounds
$\norm{\bar k_i^\delta(s)}\leq M_0$ and
$\norm{m_i h_i(s)}\leq M_0$ follow from normalization.
Moreover,
\begin{align}    
 \norm{\bar k_i^\delta(s)-m_i h_i(s)}
 &\leq\int_0^1\norm{h_i(s+\delta v)-h_i(s)}|\psi_i(v)|\,{\rm d}v\nonumber \\
 &\leq\beta_\psi M_1\delta. \label{k_bound}
\end{align}
The inequality
$\norm{aa^\top-bb^\top}
\leq(\norm{a}+\norm{b})\norm{a-b}$
and the bounds above yield
$\norm{
\bar k_i^\delta(s)\bar k_i^\delta(s)^\top
-\gamma_i h_i(s)h_i(s)^\top
}\leq 2\beta_\psi M_0M_1\delta$.
Consequently,
\[
\begin{aligned}
&\norm{\bar G_{\delta,i}-G_{0,i}}\\
&\leq
\int_0^{T-\delta}q_{i,\delta}(s)
\norm{
\bar k_i^\delta(s)\bar k_i^\delta(s)^\top
-\gamma_i h_i(s)h_i(s)^\top
}\,{\rm d}s\\
&\quad+
\gamma_i M_0^2
\int_0^T
|\tilde q_{i,\delta}(s)-q_{i,0}(s)|
\,{\rm d}s\\
&\leq
2\beta_\psi M_0M_1\delta
+M_0^2\varepsilon_\delta,
\end{aligned}
\]
where we used $\gamma_i\leq1$ and
\eqref{eq:densityerror}.
Taking convex combinations over
$p\in\Delta_m$ proves \eqref{eq:uniform}.

We next establish existence of optimizers and
uniform positive definiteness of their Gramians.
Since $\|\bar k_i^\delta(s)\|\leq M_0$, we have
$0\preceq\Gbar_\delta(p)\preceq M_0^2I_r$
for all $p\in\Pset$ and $0<\delta<T$.
On this bounded set, both objectives diverge
to $+\infty$ as the smallest eigenvalue
approaches zero.
Since $\Gbar_\delta(p^\circ)\to
\Gbar_0(p^\circ)\succ0$,
there exist $K>0$ and $\delta_0>0$ such that
$\Gbar_\delta(p^\circ)\succ0$ and
$\phi_\nu(\Gbar_\delta(p^\circ))\leq K$
for all $0<\delta<\delta_0$.
Extending the objective by $+\infty$ at
singular Gramians yields a lower semicontinuous
function on the compact set $\Pset$,
so an optimizer $p_{\nu,\delta}$ exists
\cite[Corollary~1.10]{rockafellar1998variational}.
By optimality of $p_{\nu,\delta}$,
$\phi_\nu(\Gbar_\delta(p_{\nu,\delta}))\leq K$.
For $\nu=V$, since
$\Gbar_\delta(p_{V,\delta})\preceq M_0^2I_r$,
$
e^{-K}
\leq\det\Gbar_\delta(p_{V,\delta})\leq\lambda_{\min}(\Gbar_\delta(p_{V,\delta}))
M_0^{2(r-1)}$,
which gives
$\lambda_{\min}(\Gbar_\delta(p_{V,\delta}))
\geq e^{-K}/M_0^{2(r-1)}$.
For $\nu=A$, since $Q\succ0$,
$K
\geq\lambda_{\min}(Q)
\tr(\Gbar_\delta(p_{A,\delta})^{-1})\geq\frac{\lambda_{\min}(Q)}
{\lambda_{\min}(\Gbar_\delta(p_{A,\delta}))}$.
Thus,
\begin{align}
    \Gbar_\delta(p_{\nu,\delta})\succeq c_\nu I_r, \label{key_bound}
\end{align}
where $c_V=e^{-K}/M_0^{2(r-1)}$ and
$c_A=\lambda_{\min}(Q)/K$ are positive
constants independent of $\delta$.

Let $\delta_j\downarrow0$ be arbitrary.
By compactness of $\Pset$,
the sequence $\{p_{\nu,\delta_j}\}$
has a convergent subsequence.
Passing to this subsequence, we may assume
$p_{\nu,\delta_j}\to p_*\in\Pset$.
The uniform convergence in \eqref{eq:uniform}
and \eqref{key_bound} imply
$\Gbar_0(p_*)\succeq c_\nu I_r\succ0$.
Moreover, for any $p\in\Dset_0$,
we have $\Gbar_{\delta_j}(p)\succ0$
for all sufficiently large $j$.
Hence, by optimality and continuity,
$\phi_\nu(\Gbar_0(p_*))
 =\lim_{j\to\infty}
 \phi_\nu(\Gbar_{\delta_j}(p_{\nu,\delta_j}))
 \leq\lim_{j\to\infty}
 \phi_\nu(\Gbar_{\delta_j}(p))
 =\phi_\nu(\Gbar_0(p))$.
Thus, every accumulation point minimizes
the limiting objective over $\Dset_0$.
If the limiting optimizer is unique,
every choice of $p_{\nu,\delta}$ converges
to it as $\delta\downarrow0$.
\end{IEEEproof}

\begin{remark}
The $L^1$ convergence assumption
\eqref{eq:densityerror} can be satisfied by
truncating and renormalizing a prescribed
starting-time distribution.
Specifically, given probability densities
$q_{i,0}$ on $[0,T]$, define
$q_{i,\delta}(s)
=\frac{q_{i,0}(s)}
{\int_0^{T-\delta}q_{i,0}(t)\,{\rm d}t}$
$(0\leq s\leq T-\delta)$.
In particular, if the densities $q_{i,0}$ are bounded by
a common constant $M>0$, then
$\varepsilon_\delta
=2\max_{1\le i\le m}\int_{T-\delta}^T q_{i,0}(s)\,ds
\le 2M\delta$.
Hence, $\varepsilon_\delta=O(\delta)$ as $\delta\downarrow0$.
\end{remark}

\subsection{Recovery of the Original CAS}

We now consider uniformly distributed starting times
to clarify how prescribed waveforms affect the
short-duration limit.
For $q_{i,\delta}=1/(T-\delta)$ on $[0,T-\delta]$
(and hence $q_{i,0}=1/T$ on $[0,T]$),
we have $\varepsilon_\delta=2\delta/T$ and
\begin{equation}
 \Gbar_0(p)=\frac1T\sum_{i=1}^m
 p_i\gamma_i W_i(T),
 \label{eq:waveweighted}
\end{equation}
where $W_i(T)$ is defined by \eqref{eq:cas-gramian}.
Thus, the limiting Gramian is a
waveform-efficiency-weighted version of the
original CAS Gramian.
If all waveform efficiencies $\gamma_i$ have the same positive value
$\gamma$, then
$\Gbar_0(p)=(\gamma/T)W(p,T)$,
so the optimal allocations coincide with those
of the corresponding original CAS problem~\eqref{eq:cas-objectives}.

The following corollary establishes convergence
to the original CAS under equal positive
waveform efficiencies.
The proof follows from Theorem~\ref{thm:limit}

\begin{corollary}
\label{cor:casrecovery}
Suppose $q_{i,\delta}=1/(T-\delta)$,
$\gamma_i=\gamma>0$ for every $i$, and
$W(p^\circ,T)\succ0$ for some $p^\circ\in\Pset$.
Then
\begin{equation}
 \sup_{p\in\Delta_m}
 \norm{\Gbar_\delta(p)-\gamma W(p,T)/T}
 \leq
 2M_0\left(\beta_\psi M_1+M_0/T\right)\delta.
 \label{eq:uniformcas}
\end{equation}
Every accumulation point of optimal solutions
$p_{\nu,\delta}$ to \eqref{eq:scoring2} with $Q=I_r$
as $\delta\downarrow0$ is an optimizer of the
corresponding original CAS problem
in \eqref{eq:cas-objectives}.
If the corresponding original CAS optimizer is unique,
every choice of $p_{\nu,\delta}$ converges to it
as $\delta\downarrow0$.
\end{corollary}

We next show that, when all candidates share
a common waveform, an exact correspondence
holds even for a fixed $\delta>0$.
In this case, the waveform-based Gramian can
be expressed in terms of the original
controllability Gramian at horizon $T-\delta$.

\begin{proposition}
\label{prop:exact}
Suppose that $q_{i,\delta}=1/(T-\delta)$ and all candidates share a
common waveform $\psi$.
Define
$F_{\psi,\delta}
:=\delta\int_0^1
e^{A\delta(1-v)}\psi(v)\,{\rm d}v$
and the state controllability Gramian
$X(p,t):=\int_0^t
e^{A\tau}B(p)B(p)^\top e^{A^\top\tau}
\,{\rm d}\tau$.
Then
\begin{equation}
G_\delta(p)=\frac{1}{T-\delta}
CF_{\psi,\delta}X(p,T-\delta)
F_{\psi,\delta}^\top C^\top.
\label{eq:effective}
\end{equation}
If $C=I_n$, $F_{\psi,\delta}$ is nonsingular,
and $\Dset_\delta\neq\varnothing$, the waveform-based V-CAS has
exactly the same optimal allocations as
the original V-CAS at horizon $T-\delta$.
\end{proposition}

\begin{IEEEproof}
Since $F_{\psi,\delta}$ commutes with $e^{At}$,
a change of variables gives~\eqref{eq:effective}.
For $C=I_n$, taking determinants yields
$\det G_\delta(p)
=\frac{\det(F_{\psi,\delta})^2}
{(T-\delta)^n}\det X(p,T-\delta)$.
When $F_{\psi,\delta}$ is nonsingular,
the two Gramians are positive definite
for exactly the same allocations.
Moreover, their determinants differ only
by a positive factor independent of $p$.
Hence, their V-CAS optimal allocations coincide.
\end{IEEEproof}

This exact correspondence is specific to
full-state V-CAS.
For weighted AE-CAS, the effective weight
becomes
$F_{\psi,\delta}^{-1}QF_{\psi,\delta}^{-\top}$,
so the optimal allocations need not coincide.
The nonsingularity assumption of $F_{\psi,\delta}$ holds for all
sufficiently small $\delta$ if the common
waveform has nonzero signed area
$m:=\int_0^1\psi(v)\,{\rm d}v$, because
$F_{\psi,\delta}=\delta mI_n+O(\delta^2)$.

The following corollary shows that
the original weighted AE-CAS minimizes
an asymptotic upper bound on the
target-averaged finite-trial mean-square error.

\begin{corollary}
\label{cor:finite_trial_original}
Suppose $q_{i,\delta}=1/(T-\delta)$ and
$\gamma_i=\gamma>0$ for every $i$.
For any fixed $p\in\Pset$ with
$W(p,T)\succ0$, the Stage 2 optimal policy
\eqref{eq:policy} satisfies
\begin{equation}
\limsup_{\delta\downarrow0}
\E_\eta\left[
\E_p\norm{\bar D_N-\eta}^2
\right]
\leq
\frac{TM_0^2}{N}
\tr(QW(p,T)^{-1}),
\label{eq:finite_trial_original}
\end{equation}
where $N\geq1$ is the number of independent trials,
$\bar D_N$ is their empirical mean output change
defined in Theorem~\ref{thm:sampling}, and
$M_0$ is the uniform bound on the response kernels
defined in \eqref{def_M0}.
\end{corollary}
\begin{IEEEproof}
By \eqref{eq:samplingaverage} and
$G_\delta(p)=\delta^2\Gbar_\delta(p)$,
\begin{align}
\E_\eta\left[
\E_p\norm{\bar D_N-\eta}^2
\right]
\leq
\frac{1}{N}
\left(\frac{K_\delta}{\delta}\right)^2
\tr(Q\Gbar_\delta(p)^{-1}). \label{keisan}
\end{align}
The estimate \eqref{k_bound}, the continuity
of $h_i$, and $|m_i|=\sqrt{\gamma}$ imply
$K_\delta/\delta\to\sqrt{\gamma}M_0$.
Moreover, Corollary~\ref{cor:casrecovery} and $W(p,T)\succ0$
yield
$\Gbar_\delta(p)^{-1}
\to(T/\gamma)W(p,T)^{-1}$.
Combining these limits with \eqref{keisan}
yields \eqref{eq:finite_trial_original}.
\end{IEEEproof}


\section{Numerical Example}
\label{sec:examples}

We illustrate how waveform shape and duration affect optimal allocations. Consider $A=\diag(-0.3,-3)$, $C=I_2$, $B_0=[e_1\ e_2\ e_1+e_2]$, $T=1$, $\Pset=\Delta_3$, and $Q=I_2$, with uniform starting times. We compare rectangular waveforms (R), where $\psi_i(v)=1$ for all $i$, and mixed waveforms (M), where $\psi_1(v)=\psi_2(v)=1$ and $\psi_3(v)=\sqrt{3}v$. Both satisfy the unit-$L^2$ normalization, with waveform efficiencies (defined by \eqref{efficiency}) $(1,1,1)$ and $(1,1,3/4)$, respectively. Thus, Corollary~\ref{cor:casrecovery} implies that
R recovers the original CAS in the short-duration
limit, whereas M leads to a waveform-weighted
limiting problem that differs from the original
CAS problem, as given in \eqref{eq:waveweighted}.

The candidate Gramians can be evaluated analytically. Set
$\lambda_1=0.3$, $\lambda_2=3$, and
$f_{ij}(\delta):=\int_0^1 e^{-\lambda_j\delta(1-v)}\psi_i(v)\,{\rm d}v$.
Then
 $\frac{[G_{\delta,i}]_{jk}}{\delta^2}
 =b_{ij}b_{ik}f_{ij}(\delta)f_{ik}(\delta)
 \frac{1-e^{-(\lambda_j+\lambda_k)(1-\delta)}}
 {(\lambda_j+\lambda_k)(1-\delta)}$,
where $b_{ij}$ is the $j$th component of $b_i$.
Table~\ref{tab:allocation} reports the optimal allocations for
both criteria. The rows $\delta=0$ represent the limiting
optimization problems rather than physical zero-duration
interventions.

\begin{table}[t]
\caption{Optimal allocations for R and
M waveforms.}
\label{tab:allocation}
\centering
\footnotesize
\setlength{\tabcolsep}{4pt}
\renewcommand{\arraystretch}{1.1}
\begin{tabular}{cccccccc}
\toprule
& & \multicolumn{3}{c}{V-CAS}
    & \multicolumn{3}{c}{AE-CAS}\\
\cmidrule(lr){3-5}\cmidrule(lr){6-8}
& $\delta$
& $p_1$ & $p_2$ & $p_3$
& $p_1$ & $p_2$ & $p_3$\\
\cmidrule(r){2-2}
\cmidrule(lr){3-5}
\cmidrule(l){6-8}
& 0    & 0.210 & 0.210 & 0.580
       & 0.192 & 0.588 & 0.220\\
R & 0.05 & 0.221 & 0.221 & 0.557
         & 0.187 & 0.604 & 0.209\\
& 0.40 & 0.286 & 0.286 & 0.427
       & 0.170 & 0.683 & 0.147\\
\cmidrule(r){1-2}
\cmidrule(lr){3-5}
\cmidrule(l){6-8}
& 0    & 0.305 & 0.305 & 0.390
       & 0.232 & 0.607 & 0.161\\
M & 0.05 & 0.308 & 0.295 & 0.397
         & 0.224 & 0.616 & 0.160\\
& 0.40 & 0.354 & 0.234 & 0.412
       & 0.195 & 0.661 & 0.144\\
\bottomrule
\end{tabular}
\end{table}

For R, both V-CAS and AE-CAS allocations approach the original CAS solutions as $\delta\downarrow0$. For M, the lower efficiency of the third waveform changes both limiting allocations; for example, the limiting V-CAS weight of the third candidate decreases from $0.580$ for R to $0.390$ for M. The variation across positive durations further illustrates how waveform shape and duration affect optimal intervention probabilities.

\section{Conclusion}
\label{sec:conclusion}

This paper provided an operational interpretation
of CAS through randomized finite-duration
interventions with prescribed waveforms.
We derived the minimum expected-energy policy
and the resulting allocation criteria,
and related weighted AE-CAS to finite-trial
mean-output accuracy.
We established short-duration convergence
to waveform-weighted CAS, recovering the
original CAS under equal positive efficiencies.
For common waveforms, we also established an
exact finite-duration correspondence with
full-state V-CAS.
Under uniform timing and equal positive waveform efficiencies, the original CAS provides an allocation that is insensitive to detailed waveform shapes in the short-duration limit. This property may facilitate preliminary allocation design in neural perturbation experiments before stimulation waveforms are finalized.
Future work includes noisy response estimation
and amplitude-constrained interventions.

\section*{Acknowledgment}
This work was supported by JST PRESTO, Japan, Grant Number JPMJPR25K4. 
The author used OpenAI's ChatGPT to assist with manuscript preparation. Responsibility for the final technical content and presentation rests with the author.

\bibliographystyle{IEEEtran}
\bibliography{reference}
\end{document}